\documentclass[10pt,a4paper,notitlepage, oneside]{article}
\usepackage[utf8x]{inputenc} 
\usepackage{ucs} 
\usepackage[english]{babel}
\usepackage{amsmath}
\usepackage{amsthm}
\usepackage{amsfonts}
\usepackage{cite}
\usepackage{mathtools}
\usepackage{amssymb} 
\usepackage{mathtools}
\usepackage{xfrac} 
\usepackage[percent]{overpic}
\usepackage{hyperref}
\hypersetup{
	colorlinks=true,  
	linkcolor=black,    
	urlcolor=red      
}

\usepackage{enumerate} 

\usepackage{tikz}
\usetikzlibrary{backgrounds}
\usetikzlibrary{decorations.pathmorphing,shapes}
\usetikzlibrary{arrows.meta}

\usetikzlibrary{arrows,positioning,intersections}
\usetikzlibrary{backgrounds}
\usetikzlibrary{trees}

\tikzstyle{otherbox} = [draw, thin, align=center]
\tikzstyle{minorbox} = [draw, thin, align=center, fill=gray!20]
\tikzstyle{dedge}=[thick, densely dotted]
\tikzstyle{novertex}=[rectangle]
\tikzstyle{hvertex}=[circle,inner sep=0.cm, minimum size=1mm, fill=white, draw=black]
\tikzstyle{hedge}=[ thick]
\tikzstyle{harc}=[ultra thick, ->]
\tikzstyle{point}=[draw,circle,inner sep=0.cm, minimum size=1mm, fill=black]
\tikzstyle{face}=[color=auchblau,fill=hellblau,thick]
\tikzstyle{nface}=[color=hellblau,fill=hellblau,thick] 

\colorlet{auchblau}{blue!60!white}

\colorlet{hellblau}{blue!20!white}
\colorlet{hellrot}{red!40!white}
\colorlet{hellgrau}{black!30!white}

\usepackage[textsize = tiny]{todonotes}  

\newtheorem{definition}{Definition}

\newtheorem{theorem}[definition]{Theorem}

\title{Hajós' cycle conjecture for small graphs}
\author{Irene Heinrich, Marco Natale and Manuel Streicher}
\date{}
\begin{document}

\maketitle 
	
\begin{abstract}
	Hajós' conjecture states that an Eulerian graph of order $n$ can be decomposed into at most $\lfloor \sfrac{(n-1)}{2} \rfloor$ edge-disjoint cycles.
	We describe preprocessing steps, heuristics and integer programming techniques that enable us to verify Hajós' conjecture for all Eulerian graphs with up to twelve nodes.
\end{abstract}
	
\section{Introduction}
	In 1968 Hajós conjectured that the edge set of an Eulerian graph of order $n$ can be decomposed into at most $\lfloor \sfrac{(n-1)}{2} \rfloor$ edge-disjoint cycles (cf.\ \cite{Lovasz68}).
	Granville and Moisiadis \cite{Gran87} showed that the conjecture holds true for graphs of maximum degree $4$. Furthermore, Seyffarth \cite{seyffarth1992hajos} verified the conjecture for planar graphs and, even more general, Fan and Xu \cite{Fan_Xu02} proved it for projective planar graphs.
	Moreover, Fuchs, Gellert and Heinrich \cite{Fuchs_Gellert_Heinrich} showed that Hajós' conjecture holds true for all graphs of pathwidth at most $6$.
	
	The main contribution of this paper is to verify Hajós' conjecture for small graphs:
	\begin{theorem}
		Every graph $G$ of order at most $12$ fulfils Hajós' conjecture.
	\end{theorem}
	
\begin{figure}%
\begin{tikzpicture}[node distance=2cm]
	\node (graph) {
	\begin{tikzpicture}[anchor=center, node distance=0.5cm]
	\node[minimum size=0.01cm, circle, draw] (1) {};
	\node[circle, draw, below of=1, xshift=-0.25cm] (2) {};
	\node[circle, draw, below of=1, xshift=0.25cm] (3) {};
	\node[circle, draw, below of=2, xshift=-0.25cm] (4) {};
	\node[circle, draw, below of=2, xshift=0.25cm] (5) {};
	\node[circle, draw, below of=3, xshift=0.25cm] (6) {};
	\node[left of=1, xshift=-.4cm, yshift=.35cm] (label) {Input $G$};
	\node[right of=1, xshift=1cm, yshift=.35cm] (unsichtbar) {};
	\path[-] (1) edge (2);
	\path[-] (1) edge (3);
	\path[-] (2) edge (3);
	\path[-] (2) edge (4);
	\path[-] (2) edge (5);
	\path[-] (3) edge (5);
	\path[-] (3) edge (6);
	\path[-] (4) edge (5);
	\path[-] (5) edge (6);
	\end{tikzpicture}
	};
	\node[minorbox, below of=graph, yshift=-.2cm, minimum width=9cm] (minC) {Properties of minimum counterexamples\\Fulfils all conditions of Theorem \ref{thm: minCounterex}?};
	\node[above=0mm of minC, xshift=3.7cm] {\footnotesize{Section \ref{sec: counterex}}};
	\node[minorbox, below of=minC, yshift=-.2cm, minimum width=9cm] (heu) {Random heuristics\\Found desired cycle decomposition?};
	\node[above=0mm of heu, xshift=3.7cm] {\footnotesize{Section \ref{sec:h}}};
	\node[minorbox, right of=minC, xshift=4cm, yshift=-3.45cm, rotate = 270, minimum width = 7.8cm] (out1) {Output: No minimum counterexample};
	\node[minorbox,rectangle split,rectangle split horizontal, rectangle split parts=2,below of=heu, yshift = -.5cm] (ip) {%
		\nodepart[text width=4.23cm, minimum width = 4.5cm]{one}IP-solver.\\Found feasible solution?
		\nodepart[text width=4.23cm, minimum width = 4.5cm]{two}%
		\emph{Random long cycle}.\\ Found desired cycle decomposition?};
	\node[above=0mm of ip, xshift=3.7cm] {\footnotesize{Section \ref{sec:h}}};
	\node[above=0mm of ip, xshift=-3.4cm] {\footnotesize Sections \ref{sec:ip}, \ref{sec: generalIP}};
	\node[minorbox, below=1.7cm of ip.one south, minimum width = 4.5cm] (outipno) {Output: \\Found counterexample $G$};
	\node[below of = heu, yshift = .8cm] (hilfMitte) {};
	\node[below of=heu, yshift=.55cm] (sP) {\footnotesize start parallel};
	\node[coordinate, right of = hilfMitte, xshift=0.24cm] (rechtsMitte) {};
	\node[coordinate, left of = hilfMitte, xshift=-0.24cm] (linksMitte) {};
	\node[coordinate, left=0.8cm of ip.two south] (eckeLinksOben) {};
	\node[coordinate, right=0.8cm of ip.two south] (eckeRechtsOben) {};
	\node[coordinate, below of=eckeLinksOben, yshift=1.2cm] (eckeLinksUnten) {};
	\node[coordinate, below of=eckeRechtsOben, yshift=1.2cm] (eckeRechtsUnten) {};
	\node[coordinate, right=1cm of ip.one south] (lasthelp) {};
	\node[coordinate, above=3.45cm of out1.270] (outhilf1) {};
	\node[coordinate, above=1.25cm of out1.270] (outhilf2) {};
	\node[coordinate, below=1.25cm of out1.270] (outhilf3) {};
	\node[coordinate, below=3.25cm of out1.270] (outhilf4) {};

	\path[->, thick](graph) edge (minC);
	\draw[->, thick]      (minC) -- node[fill = white] {\footnotesize Yes} (heu);
	\path[->, thick] (minC.0) edge node[fill = white] {\footnotesize No} (outhilf1);
	\path[->, thick] (heu.0) edge node[fill = white] {\footnotesize Yes} (outhilf2);
	\path[->, thick] (ip.0) edge node[fill = white] {\footnotesize Yes}  (outhilf3);
	\draw[->, thick] (heu)  |-   (rechtsMitte) -- (ip.two north);
	\draw[->, thick] (heu)  |-  node[fill = white, yshift = .35cm]{\footnotesize No} (linksMitte) -- (ip.one north);
	\draw[->, thick] (ip.one south) -- node[fill = white]{\footnotesize No} (outipno);
	\draw[->, thick] (eckeLinksOben) -- (eckeLinksUnten) -- node[fill = white]{\footnotesize No} (eckeRechtsUnten) -- (eckeRechtsOben);
	\draw[->, thick] (lasthelp) |- node[fill = white, yshift = 0.65cm] {\footnotesize Yes} (outhilf4) {};

\end{tikzpicture}
\caption{Algorithmic scheme of the verification procedure}%
\label{fig:schematicOverview}%
\end{figure}
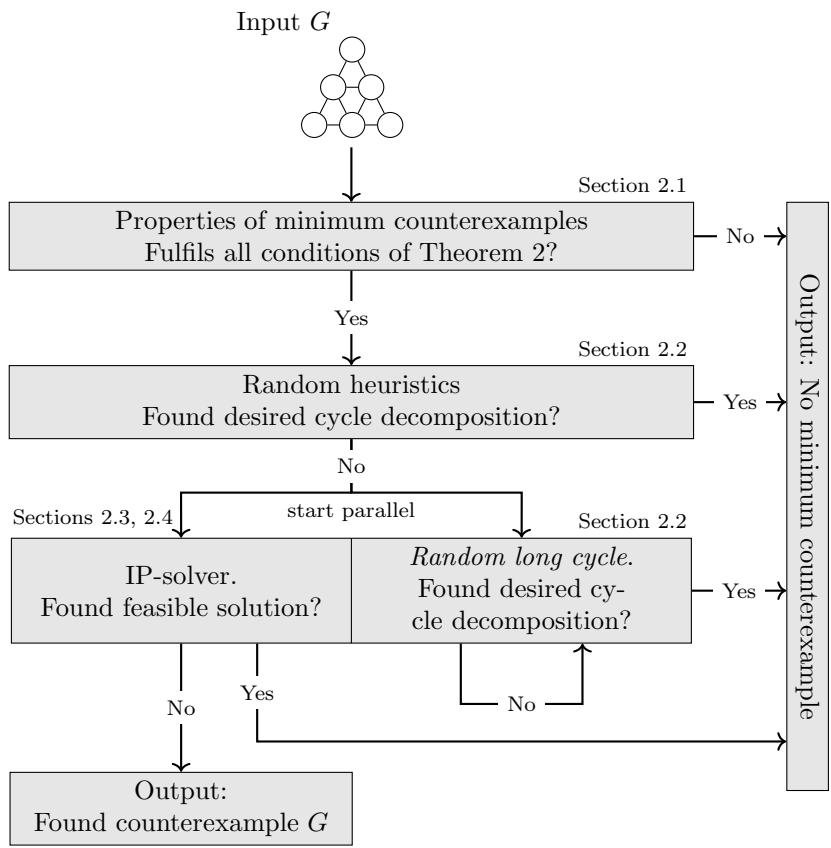

	The basic idea behind the verification procedure is as simple as it gets: We run through all graphs of order at most $12$ and test if they fulfil Hajós' conjecture.
	Nevertheless, the large number of non-isomorphic Eulerian graphs (there are $87723296$ such graphs of order $12$, cf.\ \cite{NG01}) makes this task harder than it seems.
	We use two different approaches for the tests -- minimum counterexamples and explicit computation. Figure \ref{fig:schematicOverview} illustrates the verification procedure.
	
	In Chapter \ref{sec: counterex}, essential properties of a minimum counterexample are discussed.
	Chapter \ref{sec:h} gives an insight to various heuristics which explicitly compute cycle decompositions.
	Another way of explicit calculation via integer programming techniques is then presented in Chapter \ref{sec:ip} and Appendix \ref{appendix}.
	Finally, the interaction of all the techniques and further implementation details as well as the computational results are demonstrated in Chapter \ref{sec: implementation results}.
	
\section{Methods}
\label{sec: methods}
\subsection{Minimum counterexamples to Hajós' conjecture}
\label{sec: counterex}
	The following theorem summarizes all results on minimum counterexamples which we use to verify Hajós' conjecture.
	\begin{theorem}[Minimum counterexamples to Hajós' conjecture \cite{Gran87, Fan_Xu02, Fuchs_Gellert_Heinrich}]
		\label{thm: minCounterex}
		A counterexample to Hajós' conjecture is \emph{minimum} if it is of minimum size amongst all counterexamples of minimum order.
		A minimum counterexample $G$ to Hajós' conjecture is biconnected and has the following properties:
		\begin{enumerate}[(i)]
			\item \label{itm: counterex at most 1 node <= 4} $G$ contains at most one node of degree $2$ or $4$.
			\item \label{itm: counterex adj neighb of deg 2 node} Neighbours of a degree-$2$ node in $G$ are adjacent.
			\item \label{itm: counterex neighb of deg4 node regular} The neighbourhood of a degree-$4$ node in $G$ induces a regular graph.
			\item \label{itm: counterex deg 6 node 4clique} If $v$ is a degree-$6$ node in $G$ with neighbours $\{x_1, x_2, x_3, x_4, x_5, x_6\}$ and $\{x_1, x_2, x_3, x_4\}$ is a $4$-clique in $G$, then $x_5x_6$ is an edge in $G$.
			\item \label{itm: counterex 2 adj deg6 nodes 5 common} If $G$ contains two adjacent degree-$6$-nodes $u$ and $v$ with $|N(u) \cap N(v)| = 5$. Then, $N(u) \cap N(v)$ is an independent set in $G$.
			\item \label{itm: counterex 2 adj deg6 nodes 5 common_claw} Assume that $G$ contains two adjacent degree-$6$-nodes $u$ and $v$ with $|N(u) \cap N(v)| = 5$. Then, $G-\{u,v\}$ does not contain a node with at least three neighbours in $|N(u) \cap N(v)|$.
			\item \label{itm: counterex 2 deg6 nodes 6 common} Assume that $G$ contains two adjacent degree-$6$-nodes $u$ and $v$ with $|N(u) \cap N(v)| = 4$. Let $\{x_u\} = N(u)\setminus (N(v) \cup \{v\})$ and $\{x_v\} = N(v)\setminus (N(u) \cup \{u\})$. If  $G- \{u,v\} - (N(u) \cap N(v))$ contains a path from $x_u$ to $x_v$ then, $N(u) \cap N(v)$ is an independent set in $G$.
		\end{enumerate}
	\end{theorem}
	Property \eqref{itm: counterex at most 1 node <= 4} was discovered by Fan and Xu\cite{Fan_Xu02}.
	Granville and Moisiadis \cite{Gran87} proved \eqref{itm: counterex adj neighb of deg 2 node} and \eqref{itm: counterex neighb of deg4 node regular}.
	Fuchs, Gellert and Heinrich \cite{Fuchs_Gellert_Heinrich} showed the properties \eqref{itm: counterex deg 6 node 4clique}, \eqref{itm: counterex 2 adj deg6 nodes 5 common}, \eqref{itm: counterex 2 adj deg6 nodes 5 common_claw} and \eqref{itm: counterex 2 deg6 nodes 6 common}.
	
\subsection{Heuristics}
\label{sec:h}
	We use the following heuristics to verify Hajós' conjecture. Each of which is applied recursively to the biconnected components of the given graph.
	When a cycle $C$ is found, we remove the edge set of $C$ from the given graph, add it to the cycle decomposition, and apply the heuristic again.
	
  \paragraph{Random cycle (\emph{RC})}
	Randomly walk through the graph until a cycle is closed. Return this cycle.
	\vspace{-0.3cm}
	\paragraph{Random long cycle (\emph{RLC})}
	This is a variation of \emph{RC}. Whenever a cycle is closed, save its length, but (if possible) continue walking through the graph by choosing another neighbour.
	If every current neighbour closes a cycle, return the longest cycle seen on the walk.
	\vspace{-0.3cm}
	\paragraph{Longest distance (\emph{LD})}
	Choose two nodes with maximum distance at random and find two node-disjoint paths connecting them. This can be done by a single maximum flow computation (cf.\ \cite{schrijver2002combinatorial}).
	\vspace{-0.3cm}
	\paragraph{High degree first (\emph{HDF})}
	If there is exactly one node $v$ of maximum degree, use the \emph{Random Long Cycle} heuristic starting at $v$ to find the desired cycle.
	If there are exactly two vertices of maximum degree, use the \emph{Longest Distance} heuristic to find a cycle containing both vertices.
	In all other cases, choose any node $v$ of maximum degree and use the \emph{Random Long Cycle} heuristic starting at $v$. 

\subsection{An IP-formulation for graphs with nodes of high degree}
\label{sec:ip}
	Regard an Eulerian graph $G=(V,E)$, with $n=|V|$ containing a node $\tilde{v}\in V$ with $\deg(\tilde{v}) \in \{n-2, n-1\}$.
	If $G$ fulfils Hajós' conjecture, there is a cycle decomposition of $G$ in which every cycle contains $\tilde{v}$.
	We can use this fact for a formulation as an integer linear program (IP-HD).
	For $S\subset V$ denote by $\bar{\delta}(S)$ the set of edges having at least one endpoint in $S$. Further set $I \coloneqq \{1,..,\lfloor\sfrac{(n-1)}{2}\rfloor\}$ and regard the following constraints.
	\begin{alignat}{5}
	&\hspace{-2cm}\text{(IP-HD)}&\sum_{i\in I}x_{e,i}&=1&\quad& \forall e\in E\label{ip1:1}\\
	&&\sum_{e\in\delta(v)} x_{e,i}&=2\cdot y_{v,i}&& \forall v\in V,~i\in I\label{ip1:2}\\
	&&y_{\tilde{v},i}&=1&&\forall i\in I\label{ip1:3}\\
	&&n\cdot\sum_{e\in\bar{\delta}(S)}x_{e,i}&\geq \left(n+1\right)\cdot\sum_{v\in S}y_{v,i}&&\forall i\in I,~S\subset V\setminus \{\tilde{v}\}
	&\label{ip1:4}\\
	&&x_{e,i}&\in\{0,1\} &&\forall  e\in E , i\in I\\
	&&y_{v,i}&\in\{0,1\} &&\forall v\in V, i\in I
	\end{alignat}
Variable $x_{e,i}$ has the meaning that it is equal to $1$ when edge $e$ lies on cycle $i$. Similarly, $y_{v,i}$ has the meaning that it is equal to $1$ if and only if the node $v$ lies on the $i$th cycle of the decomposition. We now show that a solution to the integer linear problem corresponds to a cycle decomposition of the desired size and vice versa.

	Let $x, y$ be a solution to the constraints above. Set $C_i\coloneqq \{e\in E \colon x_{e,i}~=~1\}$. We show that $\mathcal{C}\coloneqq \{C_i\colon i\in I\}$ corresponds to a cycle decomposition of $G$. Constraints (\ref{ip1:1}) ensure that any edge in $G$ is contained in exactly one set $C_i$. Thus, it suffices to show that the sets $C_i$ correspond to cycles in $G$. To this end fix $i\in I$ and regard equations (\ref{ip1:2}). If $y_{v,i}=1$ for some $v\in V$ we must have exactly to edges $e_1,e_2\in\delta(v)$ with $x_{e_1,i}=x_{e_2,i}=1$. On the other hand if $x_{e,i}=1$ for some $e\in E$ we must have $y_{u,i}=y_{v,i}=1$ for $(u,v)=e$. Thus, each $C_i$ is empty or corresponds to a $2$-regular graph. We now show that equations (\ref{ip1:3}) and (\ref{ip1:4}) ensure that $C_i$ corresponds to a single cycle or to the empty set. To see this, assume there exists $C'\subsetneq C_i$ corresponding to a cycle. Since $C_i$ corresponds to a $2$-regular graph there also exists $C\subsetneq C_i$ corresponding to a cycle which does not run through $\tilde{v}$. Define by $V(C)$ the nodes of this cycle. Since $\tilde{v}\notin V(C)$ equation (\ref{ip1:4}) must hold for $V(C)$. On the other hand, with (\ref{ip1:1}) we calculate
	$$n\sum_{e\in\bar{\delta}(V(C))}x_{e,i} =n\cdot\left|V(C)\right|<\left(n+1\right)\cdot\left|V(C)\right|=\left(n+1\right)\sum_{v\in V(C)}y_{v,i}.$$
	This contradicts the feasibility of $x, y$ and we may conclude that $C_i$ corresponds to a cycle or is empty.
	Thus, any feasible solution gives a cycle decomposition of $G$ into at most $\lfloor \sfrac{(n-1)}{2} \rfloor$ cycles.
	
	Now regard any cycle decomposition $\mathcal{C}$ of $G$ into at most $\lfloor \sfrac{(n-1)}{2} \rfloor$ cycles, and let $C_i$ for $1\leq i\leq \lfloor \sfrac{(n-1)}{2} \rfloor$ be the edge sets of the cycles in $\mathcal{C}$. Note that we allow $C_i$ to be the empty set. We assign values to the variables in the obvious way, $x_{e,i}=1$ if and only if $e\in C_i$ and similarly $y_{e,i}=1$ if and only if $v\in V(C_i)$.
	It is straight forward to see that equations (\ref{ip1:1}) and (\ref{ip1:2}) hold. Equations (\ref{ip1:3}) hold because $\tilde{v}$ must be contained in all cycles in the decomposition.
	Fix some $i\in I$ and $S\subset V\setminus \{\tilde{v}\}$. If $S\cap V(C_i)=\emptyset$ both sides of equation \eqref{ip1:4} equal $0$ and it is thereby fulfiled. Let now $U\coloneqq S\cap V(C_i)\neq\emptyset$. Since $\tilde{v}\in C_i$ we have $U\subsetneq V(C_i)$. Thus, we get
	$$\sum_{e\in\bar{\delta}(S)}x_{e,i}\geq 1 + \sum_{v\in S}y_{v,i}$$
	With this equation we calculate
	\begin{align*}
	n\sum_{e\in\bar{\delta}(S)}x_{e,i}&\geq n\left(1+\sum_{v\in S}y_{v,i}\right)=n + n\sum_{v\in S}y_{v,i}\\
	&\geq |S| + n\sum_{v\in S}y_{v,i}\geq (n+1)\sum_{v\in S}y_{v,i},
	\end{align*}
	which shows that all equations (\ref{ip1:4}) are fulfiled.
	
	Summing up we can test Hajós' conjecture on graphs containing a node of degree $n-1$ or $n-2$ by finding a solution to equations (\ref{ip1:1}) to (\ref{ip1:4}). 
\subsection{General IP-formulation}
\label{sec: generalIP}
Since in our experiment the general IP-formulation never managed to find a solution faster than \emph{Random long cycle}, we skip the description of the general IP here. It can be found in Appendix \ref{appendix}.
\section{Implementation and results}
\label{sec: implementation results}
In the previous section we presented the algorithms we used to verify Hajós' conjecture for all simple Eulerian graphs with twelve or less nodes.
We may restrict the verification to the class of all simple and biconnected Eulerian graphs by Theorem \ref{thm: minCounterex}.
We use the programs \texttt{geng} and \texttt{pickg} from McKay and Piperno \cite{NG01} to create all non-isomorphic, biconnected, Eulerian graphs.
All procedures introduced in Section \ref{sec: methods} are implemented in \texttt{python} version 2.7.13~\cite{python}.
In order to solve the IP-formulations, we used \texttt{Gurobi} version 7.02 \cite{gurobi}.

Let $G$ be the graph of order $n$ in the current iteration.
If $G$ does not fulfil the criteria from Theorem \ref{thm: minCounterex}, we go on to the next graph.
Otherwise, we let all heuristics presented in Section \ref{sec:h} run once and check if any heuristic finds a cycle decomposition of size less or equal than $\lfloor \sfrac{(n-1)}{2} \rfloor$.
If we still cannot confirm the conjecture for $G$ we start an IP-solver and at the same time, we let \emph{Random long cycle} run repeatedly.
If the conditions for the IP in section \ref{sec:ip} do not apply, we use the general IP-formulation (cf. Appendix \ref{appendix}).
We abort the program either when the heuristic finds a suitable cycle decomposition or when the IP-solver finishes.
Tables \ref{tab:resultsgen}, \ref{tab:resultspp} and \ref{tab:resultsh} give an overview about which procedure verified the conjecture on how many graphs.

The whole computation took roughly two days in total on a standard PC equipped with $16$~GB RAM.
However, the memory usage was always less than $2$~GB.
We stress here that the main challenge in raising the bar to graphs of higher order is not the verification for a single graph but rather the gigantic number of non-isomorphic graphs. 

\begin{table}%
\begin{center}
{\footnotesize
\begin{tabular}{c|c|c|c|c|c}
$n$&$|\mathcal{G}_n|$&Fulfilling Hajós&Theorem \ref{thm: minCounterex}&Heuristic&IP-HD\\\hline
3&1&100\%&1&0&0\\
4&1&100\%&1&0&0\\
5&3&100\%&3&0&0\\
6&7&100\%&7&0&0\\
7&30&100\%&30&0&0\\
8&162&100\%&162&0&0\\
9&1648&100\%&1645&1&2\\
10&30054&100\%&30030&11&13\\
11&1136467&100\%&1134576&1846&45\\
12&86265865&100\%&85691978&553679&20208\\
\end{tabular}
}
\caption{\footnotesize$\mathcal{G}_n$ is the set of non-isomorphic Eulerian, biconnected, simple, order $n$ graphs. Columns~$4$-$6$ show on how many graphs Hajós' Conjecture was verified using the according technique.}
\label{tab:resultsgen}
\end{center}
\end{table}

\begin{table}%
\begin{center}
{\footnotesize
\begin{tabular}{c|c|c|c|c|c|c|c|c}
$n$&Theorem \ref{thm: minCounterex}&\textbf{\eqref{itm: counterex at most 1 node <= 4}}&\textbf{\eqref{itm: counterex adj neighb of deg 2 node}}&\textbf{\eqref{itm: counterex neighb of deg4 node regular}}&\textbf{\eqref{itm: counterex deg 6 node 4clique}}&\textbf{\eqref{itm: counterex 2 adj deg6 nodes 5 common}}&\textbf{\eqref{itm: counterex 2 adj deg6 nodes 5 common_claw}}&\textbf{\eqref{itm: counterex 2 deg6 nodes 6 common}}\\\hline
3&1&1&0&0&0&0&0&0\\
4&1&1&0&0&0&0&0&0\\
5&3&3&0&0&0&0&0&0\\
6&7&7&0&0&0&0&0&0\\
7&30&29&0&0&0&1&0&0\\
8&162&159&1&0&0&0&0&2\\
9&1645&1617&1&7&8&3&0&9\\
10&30030&29442&46&282&83&22&5&150\\
11&1134576&1095272&1663&23557&5673&713&75&7623\\
12&85691978&79468073&230553&4051363&960457&54770&10047&916715\\
\end{tabular}
}
\caption{\footnotesize Detailed overview on usage of the criteria \eqref{itm: counterex at most 1 node <= 4} - \eqref{itm: counterex 2 deg6 nodes 6 common} from Theorem \ref{thm: minCounterex}.}
\label{tab:resultspp}
\end{center}
\end{table}

\begin{table}%
\begin{center}
{\footnotesize
\begin{tabular}{c|c|c|c|c|c}
$n$&Total&RC&RLC&LD&HDF\\\hline
9&1&0&1&0&0\\
10&11&0&11&0&0\\
11&1846&0&1700&11&135\\
12&553679&4&534851&687&18137\\
\end{tabular}
\caption{\footnotesize Detailed overview on usage of the heuristics from section \ref{sec:h}.}
\label{tab:resultsh}
}
\end{center}
\end{table}
\newpage

\appendix
\section{The IP-formulation}
\label{appendix}
We present here the IP-Formulation used to solve the general problem. In order to avoid the introduction of dozens of auxiliary variables and constraints we made use of so called \emph{general constraints}, which were introduced to \texttt{Gurobi} with version 7.0 (cf.~\cite{gurobi}). Among other things, they make it possible to use conjunctions, disjunctions and conditional expressions (see conditions \eqref{eq:ip2.4}-\eqref{eq:ip2.6}).
Let $G = (V, E)$ be a given Eulerian graph of order $n$.
Set $I \coloneqq \{1, 2, \dots, \lfloor\sfrac{(n-1)}{2}\rfloor\}$.
The set of feasible cycle decompositions can be described as follows: 
\begin{alignat}{3}
&\hspace{-2cm}(\text{IP-Gen})&\sum_{i\in I} x_{e,i} 			 	  &= 1	 			&\quad&\forall e\in E \label{eq:ip2.1} \\
&&\sum_{e\in \delta(v)} x_{e,i} 	&= 2y_{v,i}  	&\quad& \forall v\in V,\, i\in I \label{eq:ip2.2}\\
&&\bigvee_{v\in S} y_{v,i} &= \beta_{S,i} 	&\quad& \forall S\subset V,\, i\in I \label{eq:ip2.3} \\
&&\bigvee_{v\in V\setminus S} y_{v,i} &= \gamma_{S,i} 	&\quad&\forall S\subset V,\, i\in I \label{eq:ip2.4}\\
&&\beta_{S,i} \wedge \gamma_{S,i} &= z_{S,i} &\quad& \forall S\subset V,\, i\in I \label{eq:ip2.5}\\
&&\text{If } z_{S,i} = 1 \text{ then } \sum_{e\in \delta(S)} x_{e,i} &\geq 2 &\quad& \forall S\subset V,\, i\in I \label{eq:ip2.6}\\
&&x_{e,i} &\in \{0,1\} & &  \forall e\in E,\,i\in I \\
&&y_{v,i} &\in \{0,1\} & &  \forall v\in V,\, i\in I \\
&&\beta_{S,i},\, \gamma_{S,i},\, z_{S,i} &\in \{0,1\} && \forall S\subset V,\, i\in I.
\end{alignat}

As in Section~\ref{sec:ip} the variables $x_{e,i}$ and $y_{v,i}$ are set to one if and only if the edge $e$ and the vertex $v$ lie on the cycle corresponding to $C_i$, respectively.
Again, for a feasible solution $x,y,z,\beta,\gamma$ let $C_i\coloneqq \{e\in E\colon x_{e,i}=1 \}$. Constraints \eqref{eq:ip2.1} and (\ref{eq:ip2.2}) are the same as (\ref{ip1:1}) and (\ref{ip1:2}), and ensure that each edge is contained in exactly one $C_i$ and that the graph corresponding to $C_i$ is 2-regular. 
The constraints \eqref{eq:ip2.3} to (\ref{eq:ip2.6}) guarantee that $C_i$ indeed corresponds to a cycle. To this end, suppose that the graph corresponding to $C_i$ consists of at least two distinct cycles, say $C_i'$ and $C_i''$. 
Then we have $y_{v,i}=1$ for each vertex $v\in S\coloneqq V(C_i')$ and $\beta_{S,i}=1$. Moreover, it holds that $y_{v,i}=1$ for each $v\in V\setminus S$ and therefore $\gamma_{S,i}=1$. 
Consequently, from (\ref{eq:ip2.5}) it follows that $z_{S,i}=1$. 
Constraint (\ref{eq:ip2.6}) now yields that there exist at least two outgoing edges in $S=C_i'$, contradicting the fact that $C_i'$ corresponds to a cycle. \\
If $\mathcal{C}=\{C_i \colon i\in I\}$ corresponds to a cycle decomposition as in Section~\ref{sec:ip} we set $x_{e,i}$ and $y_{v,i} $ to one if and only if $e\in C_i$ and $v\in V(C_i)$, respectively. Again, constraint (\ref{eq:ip2.1}) and (\ref{eq:ip2.2}) are fulfilled. Let now $S\subset V$. We define $\beta_{S,i}$, $\gamma_{S,i}$ and $z_{S,i}$ as in (\ref{eq:ip2.3}), (\ref{eq:ip2.4}) and (\ref{eq:ip2.5}), respectively. If $S\cap V(C_i)\neq \emptyset$ and $(V\setminus S) \cap V(C_i) \neq \emptyset$ for some $C_i \in \mathcal{C}$ then both $\beta_{S,i}$ and $\gamma_{S,i}$ are equal to 1. Thus,  their conjunction is 1. Since $S$ as well as $V\setminus S$ contain a vertex of the cycle corresponding to $C_i$, the cut $(S,V\setminus S)$ contains at least two edges of $C_i$. Hence, constraint (\ref{eq:ip2.6}) is fulfilled. On the other hand if $S$ or $V\setminus S$ do not have any vertex in common with the cycle corresponding to $C_i$ then $z_{S,i}=0$ and (\ref{eq:ip2.6}) is trivially fulfilled.
\bibliographystyle{abbrv}

\bibliography{cycle-decomp}

\end{document}